\documentclass[11pt]{article}
\usepackage[english]{babel}
\usepackage{amsmath}
\usepackage{amsthm}
\usepackage{amsfonts}
\usepackage{amssymb}
\usepackage{enumerate}
\usepackage{xspace}
\usepackage{euscript}
\usepackage{graphicx}
\usepackage{amscd}
\usepackage{epsfig}
\usepackage{tabularx}
\usepackage{url}
\usepackage{multimedia}

\title{Recursive integral method for transmission eigenvalues}

\author{Ruihao Huang \thanks{Department of
Mathematical Sciences, Michigan Technological University, Houghton, MI 49931 ({\tt ruihaoh@mtu.edu}).}
\and Allan A. Struthers \thanks{Department of
Mathematical Sciences, Michigan Technological University, Houghton, MI 49931 ({\tt struther@mtu.edu }).}
\and Jiguang Sun \thanks{Department of
Mathematical Sciences, Michigan Technological University, Houghton, MI 49931 ({\tt jiguangs@mtu.edu}).}
\and Ruming Zhang \thanks{Department of
Mathematical Sciences, Michigan Technological University, Houghton, MI 49931 ({\tt rumingz@mtu.edu}).}}
\date{}

\begin{document}
\maketitle
\baselineskip 3.0ex
\begin{abstract}
Recently, a new eigenvalue problem, called the 
transmission eigenvalue problem,
has attracted many researchers.
The problem arose in inverse scattering 
theory for inhomogeneous media and has 
important applications in 
a variety of inverse problems for 
target identification and nondestructive testing. 
The problem is numerically challenging 
because it is non-selfadjoint and
nonlinear. 
In this paper, we propose a recursive integral method for 
computing transmission eigenvalues from a finite 
element discretization of the continuous problem. 
The method, which overcomes some difficulties
of existing methods, is based on 
eigenprojectors of compact operators. 
It is self-correcting, can separate 
nearby eigenvalues, and does not require 
an initial approximation based on some 
a priori spectral information.
These features make the method well suited 
for the transmission eigenvalue problem whose spectrum 
is complicated. Numerical examples show that the method 
is effective and robust.
\end{abstract}

\section{Introduction}
The transmission eigenvalue problem (TE)
\cite{ColtonKress2013, CakoniEtal2010IP, Sun2011IP,
CakoniGintidesHaddar2010SIAMMA}
has important applications in inverse
scattering theory for inhomogeneous media.
The problem is non-selfadjoint and not covered 
by standard partial
differential equation theory. 
Transmission eigenvalues have 
received significant attention in 
a variety of inverse problems for 
target identification and nondestructive testing
since they provide information concerning 
physical properties of the target.

Since 2010 significant effort has been 
focused on developing effective
numerical methods for transmission eigenvalues
\cite{ColtonMonkSun2010IP, Sun2011SIAMNA, 
JiSunTurner2012ACMTOM, WuCheng2013JSC, 
SunXu2013IP, AnShen2013JSC,
Kleefeld2013IP,
CakoniMonkSun2014CMAM, LiEtal2014JSC}. 
The first numerical treatment appeared in
\cite{ColtonMonkSun2010IP}, 
where three finite element methods
were proposed. A mixed method (without a
convergence proof) 
was developed in \cite{JiSunTurner2012ACMTOM}.
An and Shen \cite{AnShen2013JSC} 
proposed an efficient spectral-element based
numerical method for transmission 
eigenvalues of two-dimensional, radially-stratified media. 
The first method supported by a  
rigorous convergence analysis was introduced in
\cite{Sun2011SIAMNA}.
In this article transmission eigenvalues are computed 
as roots of a nonlinear function whose values are 
eigenvalues of a related positive definite 
fourth order problem. This method has two drawbacks 
1) only real transmission eigenvalues can
be obtained, and 2) many fourth order 
eigenvalue problems need to be solved. 
In \cite{Kleefeld2013IP} (see also \cite{Beyn2012LAA})
surface integral and contour integral based methods 
are used to compute both real and complex 
transmission eigenvalues in the special case 
when the index of refraction is constant. 
Recently, Cakoni et.al. \cite{CakoniMonkSun2014CMAM} 
reformulated the problem and proved convergence
(based on Osborn's compact operator theory 
\cite{Osborn1975MC}) of a mixed finite 
element method. Li et.al. 
\cite{LiEtal2014JSC}
developed a finite element method based 
on writing the TE as a quadratic eigenvalue problem.
Some non-traditional methods, including
the linear sampling method in the inverse
scattering theory \cite{Sun2012IP}
and the inside-out duality 
\cite{LechleiterRennoch2015SIAMMA}, were proposed to search 
for eigenvalues using scattering data. 
However, these methods seem to be computationally 
prohibitive since they rely on solving 
tremendous numbers direct problems.
Other methods 
\cite{GintidesPallikarakis2013IP, 
CossonniereHaddar2013JIEA, JiSun2013JCP,
JiSunXie2014JSC} and the related 
source 
problem \cite{HsiaoEtal2011CAM, WuCheng2013JSC}
have been discussed in the literature. 

In general, developing effective finite element 
methods for transmission eigenvalues is challenging
because it is a quadratic, typically, 
degenerate, non-selfadjoint eigenvalue
problem for a system of two second order partial
differential equations and despite 
some qualitative estimates the spectrum 
is largely unknown. In most cases, the continuous 
problem is degenerate with an infinite dimensional
eigenspace associated with a zero eigenvalue.
The system can be reduced to
a single fourth order problem however conforming 
finite elements for such problems
{\it e.g.} Argyris are expensive.
Straightforward finite element discretizations
generate computationally challenging large sparse 
non-Hermitian matrix eigenvalue problems. Traditional 
methods such as shift and invert Arnoldi are 
handicapped by the lack of a priori eigenvalue 
estimates.  To summarize, finite element 
discretizations of transmission eigenvalue problems
generate large, sparse, typically highly degenerate,
non-Hermitian matrix eigenvalue problems with little 
a priori spectral information beyond the likelihood 
of a relatively high-dimensional nullspace.

These characteristics suggest that most existing 
eigenvalue solver are unsuitable for transmission 
eigenvalues.  Recently integral based 
methods \cite{SakuraiSugiura2003CAM, 
Polizzi2009PRB}
related to the earlier work
\cite{Goedecker1999RMP}
and originally developed for electronic structure 
calculations become popular. These
methods are based on eigenprojections 
\cite{Kato1966}
provided by contour integrals of
the resolvent 
\cite{AustinKravajaTrefethen2014SIAMNA}.

In this paper, we propose a recursive integral method 
({\bf RIM}) to compute transmission eigenvalues 
from a continuous finite element discretization.
Regions in the complex plane are searched for 
eigenvalues using approximate eigenprojections
onto the eigenspace associated with the eigenvalues 
within the region. The approximate eigenprojections
are generated by approximating the resolvent 
contour integral around the boundary of the region 
by a quadrature on a random sample. The region is 
subdivided and subregions containing
eigenvalues are recursively subdivided 
until the eigenvalues are localized to the 
desired tolerance. 
{\bf RIM} is designed to approximate all 
eigenvalues
within a specific region without 
resolving eigenvectors. This is well
suited to the transmission eigenvalue 
problem which typically seeks only
the eigenvalues near but not at 
the origin. The degenerate
non-hermitian nature of the matrix and
the complicated unknown structure 
of the spectrum are not an issue.

{\bf RIM} is distinguished from other 
integral methods in literature by several features.
First, the method works for Hermitian and non-Hermitian 
generalized eigenvalue problems such as those from
the discretization of non-selfadjoint partial differential 
equations. Second, the recursive procedure automatically 
resolves eigenvalues near region boundaries and
minimally separated eigenvalue pairs. Third, 
the method requires only linear solves with 
no need to explicitly form a matrix inverse.

The paper is arranged as follows.
Section 2 introduces the transmission eigenvalue 
problem, the finite element discretization, and
the resulting large sparse non-Hermitian generalized
matrix eigenvalue problem. 
Section 3 introduces the recursive 
integral method {\bf RIM} to compute 
all eigenvalues within a region of the complex 
plane. 
Section 4 details various implementation details.
Section 5 contains results from a range of numerical 
examples.
Section 6 contains discussion and future work.

\section{The transmission eigenvalue problem}
\subsection{Formulation}
We introduce the transmission eigenvalue problem related to the Helmholtz equation.
Let $D \subset \mathbb R^2$ be an open bounded domain with a Lipschitz boundary $\partial D$.
Let $k$ be the wave number of the incident wave 
$u^i = e^{ikx\cdot d}$ and $n(x)$ be the index of refraction. The direct scattering problem is
to find the total field
$u(x)$ satisfying
\begin{subequations}
\begin{align}
 &\nabla \cdot  \nabla u+k^2n(x)u=0,  &\text{in } D,\\[1mm]
 &\Delta u+k^2u=0, &\text{in } \mathbb R^2 \setminus D,\\[1mm]
 &u(x)=e^{ikx\cdot d}+u^s({x}),  &\text{on } \mathbb R^2 , \\[1mm]
\label{Sommer1} & \lim_{r\to\infty}\sqrt{r} \left(\frac{\partial u^s}{\partial r}-iku^s \right)=0,
\end{align}
\end{subequations}
where $u^s$ is the scattered field, $x \in \mathbb R^2 , r = |x|$, $d\in \Omega:=\{\hat{x}\in\mathbb R^2 ; |\hat{x}|=1\}$.
The Sommerfeld radiation condition \eqref{Sommer1} is assumed to hold uniformly with respect to  $\hat{x} = x/|x|$.  

The associated transmission eigenvalue problem is to find $k$ such that there exist non-trivial solutions $w$ and $v$ satisfying
\begin{subequations}\label{ATE}
\begin{align}
\label{ATEw}&\nabla \cdot \nabla w+k^2n(x)w=0, &\text{in } D,\\[1mm]
\label{ATEv}&\Delta v+k^2v=0, &\text{in } D,\\[1mm]
\label{ATEbcD}&w -v = 0, &\text{on } \partial D, \\[1mm]
\label{ATEbcN}& \frac{\partial w}{\partial \nu} - \frac{\partial v}{\partial \nu}= 0, &\text{on } \partial D,
\end{align}
\end{subequations}
where $\nu$ the unit outward normal to $\partial D$.
The wave numbers $k$'s for which the transmission 
eigenvalue problem has non-trivial solutions 
are called transmission eigenvalues.
For existence results for 
transmission eigenvalues the reader is referred to the article
and reference list of
\cite{CakoniGintidesHaddar2010SIAMMA}.

It is clear that $k=0$ and $w = v$ a 
harmonic function in $D$ satisfies
\eqref{ATE}. So $k=0$ is a non-trivial 
transmission eigenvalue with an infinite 
dimensional eigenspace.


\subsection{A continuous finite element method}
In the following, we describe a continuous finite element 
method for \eqref{ATE} 
\cite{CakoniEtal2010IP, JiSun2013JCP}.
We use standard linear Lagrange finite element for
discretization and define
\begin{eqnarray*}
S_h&=&\mbox{ the space of continuous
piecewise linear finite elements on }D,\\
S_{h}^{0}&=&S_h\cap H_0^1(D)\\
&=&  \mbox{ the subspace of functions in } 
S_h \mbox{ with vanishing DoF  on }\partial D,\\
S_{h}^{\mathcal{B}}&=&
 \mbox{ the subspace of functions in } 
 S_h \mbox{ with vanishing DoF  in }D,
\end{eqnarray*}
where DoF stands for degrees of freedom. 

Multiplying \eqref{ATEw} by a test function $\phi$ and 
integrating by parts gives
\begin{equation}\label{ATEwWeak}
(\nabla w,  \nabla \phi) - k^2 ( n w, \phi) - \left\langle \frac{\partial w}{\partial \nu}, \phi \right\rangle = 0,
\end{equation}
where $\langle \cdot, \cdot \rangle$ denotes 
the boundary integral on $\partial D$.
Similarly, multiplying \eqref{ATEv} by a 
test function $\phi$ and integrating by parts 
gives
\begin{equation}\label{ATEvWeak}
( \nabla v,  \nabla \phi ) - k^2 ( v, \phi) - \left\langle \frac{\partial v}{\partial \nu},  \phi \right\rangle  = 0.
\end{equation}
Subtracting \eqref{ATEvWeak} from 
\eqref{ATEwWeak} and using 
the boundary condition \eqref{ATEbcN} gives
\begin{equation}\label{ATEvwWeak}
( \nabla w - \nabla v, \nabla \phi ) - 
k^2 ((nw-v), \phi ) = 0.
\end{equation}

The Dirichlet boundary condition 
\eqref{ATEbcD} 
is explicitly enforced on the discretization 
by setting
\begin{eqnarray*}
w_h&=&w_{0,h}+w_{\mathcal{B},h}
\mbox{ where }w_{0,h}\in S_h^0
\mbox{ and } w_{\mathcal{B},h}\in
S_{h}^\mathcal{B},\\
v_h&=&v_{0,h}+w_{\mathcal{B},h}
\mbox{ where }v_{0,h}\in S_h^0.
\end{eqnarray*}
Choosing the test function  
$\xi_h \in S_h^0$ 
for \eqref{ATEwWeak} gives
the weak formulation for $w_{h}$ as
\begin{equation}\label{Anablaw0h}
(\nabla(w_{0,h}+w_{\mathcal{B},h}), \nabla\xi_h )-k^2(n(w_{0,h}+w_{\mathcal{B},h}), \xi_h)=0,
\end{equation}
for all $\xi_h\in S_h^0$.
Similarly, choosing the test function 
$\eta_h\in S_h^0$ gives the weak
formulation for $v_{h}$ as
\begin{equation}\label{nablav0h}
(\nabla(v_{0,h}+w_{\mathcal{B},h}),\nabla\eta_h)-k^2((v_{0,h}+w_{\mathcal{B},h}), \eta_h)=0,
 \end{equation}
for all $\eta_h\in S_h^0$.
Finally, choosing $ \phi_h\in S_{h}^\mathcal{B}$  in \eqref{ATEvwWeak} gives
\begin{eqnarray} \nonumber
&&(\nabla (w_{0,h}+w_{\mathcal{B},h}),\nabla\phi_h)-(\nabla (v_{0,h}+w_{\mathcal{B},h}),\nabla\phi_h\,) \\
\label{Aw0h} && \qquad \qquad \qquad  -k^2\left( n(w_{0,h}+w_{\mathcal{B},h})-(v_{0,h}+w_{\mathcal{B},h}),\phi_h \right)=0.
\end{eqnarray}

Let  $\{\xi_1, \ldots, \xi_{N_h^0} \}$
be the finite element basis for $S_h^0$ and 
\[
\{\xi_1, \ldots, \xi_{N_h^0},  \xi_{N_h^0+1}, \ldots, \xi_{N_h}\}
\] 
be the basis for $S_h$. Let $N_h$, $N_h^0$, and $N_h^\mathcal{B}$ be the dimensions of $S_h$, $S_h^0$ and $S_h^\mathcal{B}$, respectively.
Clearly $\{ \xi_{N_h^0+1}, \ldots, \xi_{N_h}\}$ is 
a basis for $S_h^\mathcal{B}$ and
\[
N_h = N_h^0 + N_h^\mathcal{B}.
\]

Let $S$ be the  stiffness matrix given by 
$(S)_{j,\ell}=(\nabla\xi_j,
\nabla\xi_{\ell})$, 
$M_n$ be the mass matrix given by 
$(M_n)_{j,\ell}=(n\xi_j,\xi_{\ell})$, and
$M$ be the mass matrix given by 
$(M)_{j,\ell}=( \xi_j,\xi_{\ell})$.
Combining \eqref{Anablaw0h}, 
\eqref{nablav0h}, and \eqref{Aw0h}, gives the
generalized eigenvalue problem
\begin{equation}\label{cfemgep1}
A{\bf x}=k^2 B {\bf x},
\end{equation}
where matrices ${A}$ and ${B}$ are
\[
{A}=\left(\begin{array}{ccc}
S^{N_h^0\times N_h^0}&0&S^{N_h^0 \times N_h^\mathcal{B}} \\
0&S^{N_h^0\times N_h^0}&S^{N_h^0 \times N_h^\mathcal{B}}\\
(S^{N_h^0 \times N_h^\mathcal{B}})^T&(-S^{N_h^0 
\times N_h^\mathcal{B}})^T&S^{N_h^\mathcal{B} 
\times N_h^\mathcal{B}}-S^{N_h^B
\times N_h^\mathcal{B}}\end{array}\right),
\]
and
\[
{B}=\left(\begin{array}{ccc}
M_n^{N_h^0\times N_h^0}&0&M_n^{N_h^0\times N_h^\mathcal{B}}\\
0&M^{N_h^0\times N_h^0}&M^{N_h^0\times N_h^\mathcal{B}}\\
(M_n^{N_h^0\times N_h^\mathcal{B}})^T&-(M^{N_h^0
\times N_h^\mathcal{B}})^T&M_n^{N_h^\mathcal{B} 
\times N_h^\mathcal{B}}-M^{N_h^B
\times N_h^\mathcal{B}}\end{array}\right).
\]

${A}$ and ${B}$ are clearly not symmetric and
in general there are complex eigenvalues.
Applications are typically interested in 
determining the structure of the spectrum (including
complex conjugate pairs) near the origin. 
In practice, the primary focus is on
computing a few of the
non-trivial eigenvalues nearest the origin. 
Note for the transmission eigenvalue problem 
eigenvectors are of significantly less interest.

Arnoldi iteration based adaptive search methods for real
transmission eigenvalues were developed in 
\cite{JiSunTurner2012ACMTOM} and
\cite{MonkSun2012SIAMSC}.
However, these methods are inefficient, may fail
to converge, and are unable to compute all 
eigenvalues in general.
The main goal
of the current paper is to develop an effective 
tool to compute all the 
transmission eigenvalues (real and complex) 
of \eqref{cfemgep1} in a region of the complex plane.

\section{A recursive contour integral method}
\subsection{Continuous case}
We start by recalling some classical results in operator theory (see, e.g., \cite{Kato1966}).
Let $T: \mathcal{X} \to \mathcal{X} $ be a compact operator on a complex Banach space $\mathcal{X} $. The resolvent of $T$ is defined as
\begin{equation}\label{rhoT}
\rho(T)=\{ z \in \mathbb C: (z-T)^{-1} \text{ exists as a bounded operator on } \mathcal{X} \}.
\end{equation}
For any $z \in \rho(T)$, 
\[
R_z(T) = (z-T)^{-1}
\] 
is the resolvent of $T$ and the spectrum of $T$ is $\sigma(T)=\mathbb C
\setminus \rho(T)$.

Let $\Gamma$ be a simple closed curve on the complex plane $\mathbb C$ lying in $\rho(T)$ which contains $m$ eigenvalues of $T$: $\lambda_i, i = 1, \ldots, m$. 
The spectral projection
\[
E(T)= \frac{1}{2\pi i} \int_\Gamma R_z(T) dz.
\]
is a projection onto the space of generalized
eigenfunctions $u_i, i=1, \ldots, m$ 
associated with the eigenvalues $\lambda_i, i=1, \ldots,m$.
If a function $f$ has components in 
$u_i, i=1, \ldots, m$ then $E(T) f $ is non-zero. 
If $f$ has no components in $u_i, i=1, \ldots, m$
then $E(T)f = 0$. Thus $E(T)f$ can
be used to decide if a region contains eigenvalues of 
$T$ or not. This is the basis of {\bf RIM}.

Our goal is to compute all the 
eigenvalues of $T$ in a region $S\subset \mathbb C$.
{\bf RIM}
starts by defining $\Gamma = \partial S$,
randomly choosing several functions 
$f_j, j = 1, \ldots, J$ and approximating 
\[
I_j = E(T) f_j, \quad j=1, \ldots, J,
\]
by a suitable quadrature. Based on $I_j$  we decide if
there are eigenvalues inside $S$. If $S$ contains
eigenvalue(s), we partition $S$ into subregions 
and recursively repeat this procedure for each 
subregion. The process terminates when each
eigenvalue is isolated within a sufficiently 
small subregion. 

\begin{itemize}
\item[] {\bf RIM}$(S, \epsilon, f_j, j=1,\ldots, J)$
\item[{\bf Input:} ] search region $S$, tolerance $\epsilon$, random functions $f_j, j=1,\ldots, J$
\item[{\bf Output:} ] $\lambda$, eigenvalue(s) of $T$ in $S$
\item[1.] Approximate (using a suitable quadrature) 
the integral
		\[
		 E(T) f_j = \frac{1}{2\pi i} \int_\Gamma R_z(T) dz f_j, \quad j=1, \ldots, J, \quad \Gamma = \partial S.
		\]
\item[2.] Decide if $S$ contains eigenvalue(s):
	\begin{itemize}
		\item No. exit.
		\item Yes. compute the size $h(S)$ of $S$
			\begin{itemize}
				\item[-] if $h(S) > \epsilon $, partition $S$ into subregions $S_i, i=1, \ldots N$
						\begin{itemize}
						\item[] for $i=1$ to $N$
						\item[] $\qquad {\bf RIM}(S_i, \epsilon, f_j, j=1, \ldots, J)$
						\item[] end
						\end{itemize}
				\item[-] if $h(S) \le \epsilon$, output the eigenvalue $\lambda$ and exit
			\end{itemize}
	\end{itemize}
\end{itemize}

\subsection{Discrete case}
We specialize {\bf RIM} to potentially
non-Hermitian generalized matrix eigenvalue 
problems. The finite element discretization
of the transmission eigenvalue problem
produces such a problem as do other similar
discretizations of other PDEs. 

The matrix eigenvalue problem is
\begin{equation}\label{AxlambdaBx}
A{\bf x} = \lambda B{\bf x}
\end{equation}
where $A, B$ are $n \times n$ matrices,
$\lambda$ is a scalar, and 
${\bf x}$ is an $n \times 1$ vector. 
The resolvent is
\begin{equation}\label{eigresolvent}
R_z(A, B) = (zB-A)^{-1}B
\end{equation}
for $z$ in the resolvent set of the matrix pencil.
The projection onto the generalized 
eigenspace corresponding to eigenvalues 
enclosed by a simple closed curve $\Gamma$
is given by the Cauchy integral
\begin{equation}\label{ABpro}
E(A, B) = \frac{1}{2\pi i} \int_\Gamma (zB-A)^{-1}B dz .
\end{equation}

If the matrix pencil is non-defective then
$AX = BX\Lambda$ where $\Lambda$ is a 
diagonal matrix of eigenvalues and 
$X$ is an invertible matrix of generalized eigenvectors.
This eigenvalue decomposition shows
\[
(zB - A) X = z B X - A X = z B X-B X \Lambda = B X 
(z I - \Lambda),
\]
and gives
\[
X (zI - \Lambda)^{-1} X^{-1} = (zB - A)^{-1} B
\]
for complex $z$ not equal to any of the eigenvalues. 
Integrating the resolvent around a closed contour 
$\Gamma$ in $\mathbb C$ gives
\[
	\dfrac{1}{2 \pi i} \int_{\Gamma} R_z(A, B){d} z 
	= X \dfrac{1}{2 \pi i} \int_{\Gamma} 
	(zI - \Lambda)^{-1} {d} z X^{-1} 
	= X \Lambda_{\Gamma} X^{-1},
\]
where $\Lambda_{\Gamma}$ is $\Lambda$ with eigenvalues
inside $\Gamma$ set to $1$ and those outside $\Gamma$ set to
$0$. 

The projection of a vector ${\bf y}$ onto the 
generalized eigenspace for eigenvalues inside $\Gamma$ is
\begin{equation}\label{XLXy}
	P{\bf y}:=X \Lambda_{\Gamma} X^{-1} {\bf y} = \dfrac{1}{2 \pi i} \int_{\Gamma} R_z(A, B){\bf y} {d} z.
\end{equation}
If there no eigenvalues are inside $\Gamma$, 
then $P=0$ and  $P{\bf y} = {\bf 0}$ for all 
${\bf y} \in \mathbb C^n$.

We select a quadrature rule 
to approximate the contour integral
\[
\frac{1}{2 \pi i} \int_{\Gamma} R_z(A, B) {\bf y} {d} z \approx \frac{1}{2 \pi
i} \sum_{q=1}^N \omega_q R_{z_q}(A, B) {\bf y},
\]
where $\omega_q$ and $z_q$ are 
the quadrature weights and points, respectively.
Although an explicit computation of $R_z$ is not
possible one can approximate the projection of
${\bf y}$ by
\begin{equation}\label{Py}
P {\bf y} \approx  \sum_{q=1}^N {\bf r}_q.
\end{equation}
where ${\bf r}_q$ are the solutions
of the linear systems
\begin{equation}\label{projectionr}
(z_q B-A) {\bf r}_q = 
\frac{1}{2\pi i} \omega_q B {\bf y}, \quad q=1, \ldots, N.
\end{equation}

For robustness, we use a set of vectors 
${\bf y}_j, j=1,\ldots, J$ assembled as 
the columns of an $n \times J$ matrix $Y$.
The {\bf RIM} for generalized eigenvalue problems is as follows.
\begin{itemize}
\item[] {\bf M-RIM}$(A, B, S, \epsilon, Y)$
\item[{\bf Input:} ] matrices $A$ and $B$, search region $S$, tolerance $\epsilon$, random vectors $Y$
\item[{\bf Output:} ] generalized eigenvalue $\lambda$
\item[1.] Compute 
$P{\bf y}_j, j = 1, \ldots, J $ using \eqref{Py} on
$\partial S$.
\item[2.] Decide if $S$ contains eigenvalue(s):
	\begin{itemize}
		\item No. exit.
		\item Yes. compute the size $h(S)$ of $S$
			\begin{itemize}
				\item[-] if $h(S) > \epsilon $, partition $S$ into subregions $S_i, i=1, \ldots I$
						\begin{itemize}
						\item[] for $i=1$ to $I$
						\item[] $\qquad${\bf M-RIM}$(A, B, S_i, \epsilon, Y)$
						\item[] end
						\end{itemize}
				\item[-] if $h(S) \le \epsilon$, output the eigenvalue $\lambda$ and exit
			\end{itemize}
	\end{itemize}
\end{itemize}

\section{Implementation}
We assume the search region $S$ is a polygon in 
the complex plane $\mathbb C$
for simplicity and divide $S$ into 
subregions of simple geometry, such as
triangles and rectangles. 
Rectangles are used in the implementation.

There are several keys in the implementation of {\bf RIM}: we need a suitable quadrature rule for 
the contour integral; we need a mechanism 
to solve \eqref{projectionr}; and
we need an effective rule to decide if 
a subregion contains eigenvalues.

We use Gaussian quadrature on each rectangle 
edge. It does not appear necessary to use 
many points and we use the two point rule.

In contrast with the quadrature, an accurate linear solver seems
necessary and we use the Matlab ``$\backslash$''
command.

%
Next we discuss the rule to decide if
$S$ might contain eigenvalues 
and needs
to be subdivided. We refer to a subregion
that potentially contains at least one 
eigenvalue as admissible. Any vector
${\bf y}$ is represented in the eigenbasis
(columns of $X$) as 
${\bf y}=\sum_{i=1}^{n} a_i \textbf{x}_i $.
Assume there are $M$ eigenvalues inside 
$\Gamma$ and reorder the eigenvalues and eigenvectors
with these $M$ eigenvectors as $\textbf{x}_1 , 
\textbf{x}_2 , \dots,\textbf{x}_M$ then
\[
P{\bf y}= \dfrac{1}{2 \pi i} \int_{\Gamma} R_z(A, B) {\bf y} dz=X\Lambda_{\Gamma}X^{-1}{\bf y}=\sum_{i=1}^M a_i \textbf{x}_i .
\]
So it is reasonable to use $\|P{\bf y}\|$ to decide 
if a region contains eigenvalues. 
There are two primary concerns for the robustness
of the algorithm. We might miss eigenvalues if 
$\|P{\bf y}\|$ is small when there is an eigenvalue
within $\Gamma$.  We might continue to subdivide
a region if $\|P{\bf y}\|$ is large when there is no 
eigenvalue within $\Gamma$. In the first case 
$\|P{\bf y}\|$ could be small when there is an
eigenvalue because of quadrature/rounding errors 
and/or simply because the random components 
$a_i$ are small. Our solution is to project $P{\bf y}$ 
again with an amplifier $K$ and look at $\|P(KP{\bf y})\|$. 
In fact, one can simply choose $K = 1/\|P{\bf y}\|$.
In the second case  $\|P{\bf y}\|$ could be large 
when there is no eigenvalue inside $\Gamma$
if there are eigenvalues right outside $\Gamma$ and 
the quadrature rule or the linear solver are not 
sufficiently accurate. Fortunately,
{\bf RIM} has an interesting 
{\it self-correction} property 
that fixes such errors on subsequent iterations.  

In our implementation, we use the following rules to 
decide an admissible region:
\begin{itemize}
\item[1.] We use several random vectors ${\bf y}_j, j=1, \ldots, J$;
\item[2.] We use $\|P(KP{\bf y}_j)\|$ where $K$ is an amplifier.
\end{itemize}
Rule 1. and Rule 2. guarantee that even if the 
component of ${\bf y}$ in $X$ is small, 
the algorithm can detect it effectively
since $P(KP{\bf y})$ should be of the 
same size of $KP({\bf y})$.
If there is no eigenvalue inside $\Gamma$, $P{\bf y}$ 
can still be large due to
reasons we mentioned above.
However, another projection 
of $P(KP{\bf y})$ should significantly reduce $\|KP{\bf y}\|$.

The indicator function $\chi_S$ is the ratio
\[
\chi_S := \dfrac{ \Vert P(KP{\bf y})\Vert }{ \Vert P{\bf y} \Vert}.
\]
If there are eigenvalues inside $\Gamma$, then $\dfrac{ \Vert P(KP{\bf y})\Vert }{ \Vert P{\bf y} \Vert} = O(K) $.
On the contrary, if there is no eigenvalue inside $\Gamma$, $\dfrac{ \Vert P(KP{\bf y})\Vert }{ \Vert P{\bf y} \Vert} =o(K) $.

Here are some details in the actual implementation.
\begin{enumerate}
\item[1.] The search region $S$ is a rectangle.
\item[2.] We use $3$ random vectors ${\bf y}_j,j=1, 2, 3$.
\item[3.] The amplifier is set as $K=10$.
\item[4.] We use $2$ point Gauss quadrature rule on each edge of $S$.
\item[5.] We use Matlab "$\backslash$" to solve the linear systems.
\item[6.] We take the indicator function as
	\[
	\chi_S = \max_{j=1, 2, 3} \dfrac{\Vert P(K P{\bf y}_j)) \Vert}{ \Vert P{\bf y}_j \Vert}.
	\]
\item[7.] We use $K/10$ as the criterion, i.e., if $\chi_S > K/10$, $S$ is
admissible.
\end{enumerate}


\section{Numerical Examples}
In this section, we assume that the initial 
search region $S$ is a rectangle. 
We present examples to show the performance of {\bf RIM}.

\subsection{Transmission Eigenvalues}
We test {\bf RIM} on the generalized 
matrix eigenvalue problem for transmission 
eigenvalues using continuous finite element 
method described in Section 2. Since the 
original partial differential problem is non-selfadjoint,
the generalized matrix eigenvalue problem is non-Hermitian. 
In practice, we only need a few
eigenvalues of smallest norm. 
However, we do not have an a prior knowledge of the locations of the eigenvalues.

{\bf Example 1:} We consider a disc 
$D$ with radius $1/2$ and index of refraction 
$n(x) = 16$ where the exact transmission eigenvalues
\cite{ColtonMonkSun2010IP} are the roots of 
\[
J_1(k/2)J_0(2k) = 4 J_0(k/2)J_1(2k), \quad m = 0,
\]
and
\[
J_{m-1}(k/2) J_m(2k) = 4J_m(k/2) J_{m-1}(2k), \quad m \ge 1,
\]
where $J_m$'s are Bessel functions.

A regular mesh with with $h \approx
0.05$ is used to generate the $1018 \times 1018$
matrices $A$ and $B$ and we consider the 
preliminary search region $S=[1, 10] \times [-1, 1]$.
Since the mesh is relatively coarse
we take $\epsilon = 1.0e-3$ and 
use 3 random vectors. {\bf RIM} computes $3$ eigenvalues 
\[
\lambda_1 = 3.994, \quad \lambda_2 = 6.935, \quad \lambda_3 = 6.939
\]
which are good approximations of the exact eigenvalues given in
\cite{CakoniEtal2010IP}
\[
\lambda_1 = 3.952, \quad \lambda_2 = 6.827, \quad \lambda_3 = 6.827.
\]
Note that the values we compute are $k^2$'s and the actual values in \cite{CakoniEtal2010IP} are $k$'s.

As a second test we choose 
$S=[22, 25] \times[-8, 8]$ and 
find $4$ eigenvalues in this region
\[
\lambda_1 = 24.158+ 5.690 i, \quad \lambda_2 = 24.158- 5.690 i,\quad \lambda_3 =
25.749, \quad \lambda_4 = 25.692
\]
which approximate the exact eigenvalues
\[
\lambda_{1,2} = 23.686\pm 5.667 i, \quad \lambda_{3,4} = 24.465.
\]
Note that {\bf RIM} computes the generalized eigenvalues
to the anticipated accuracy $\epsilon$ the discrepancy
is mainly due to the fact finite element methods
approximate smaller eigenvalues better than larger eigenvalues.

\begin{figure}[hh]
\begin{center}
\begin{tabular}{lll}
\scalebox{0.42}{\includegraphics{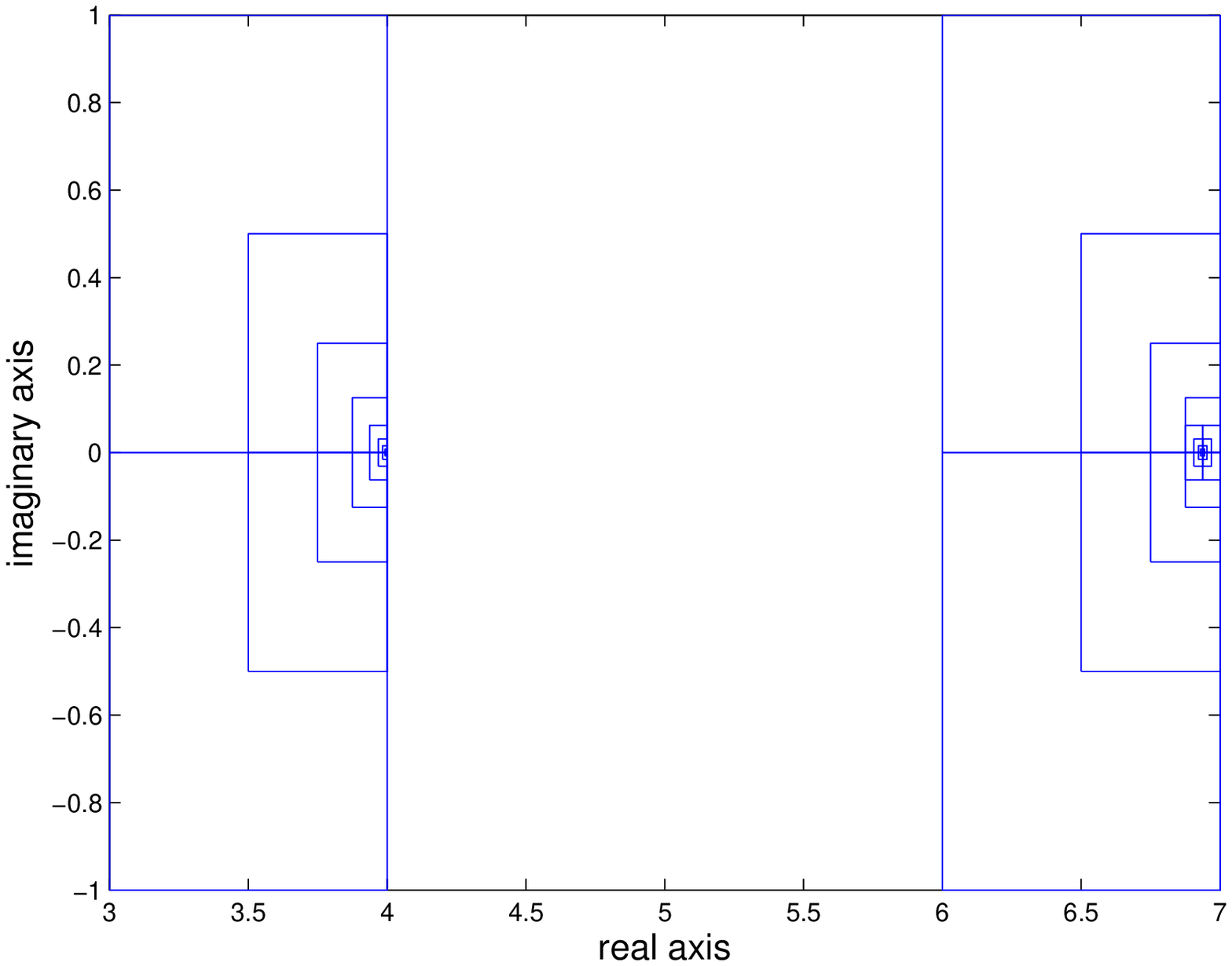}}&
\scalebox{0.42}{\includegraphics{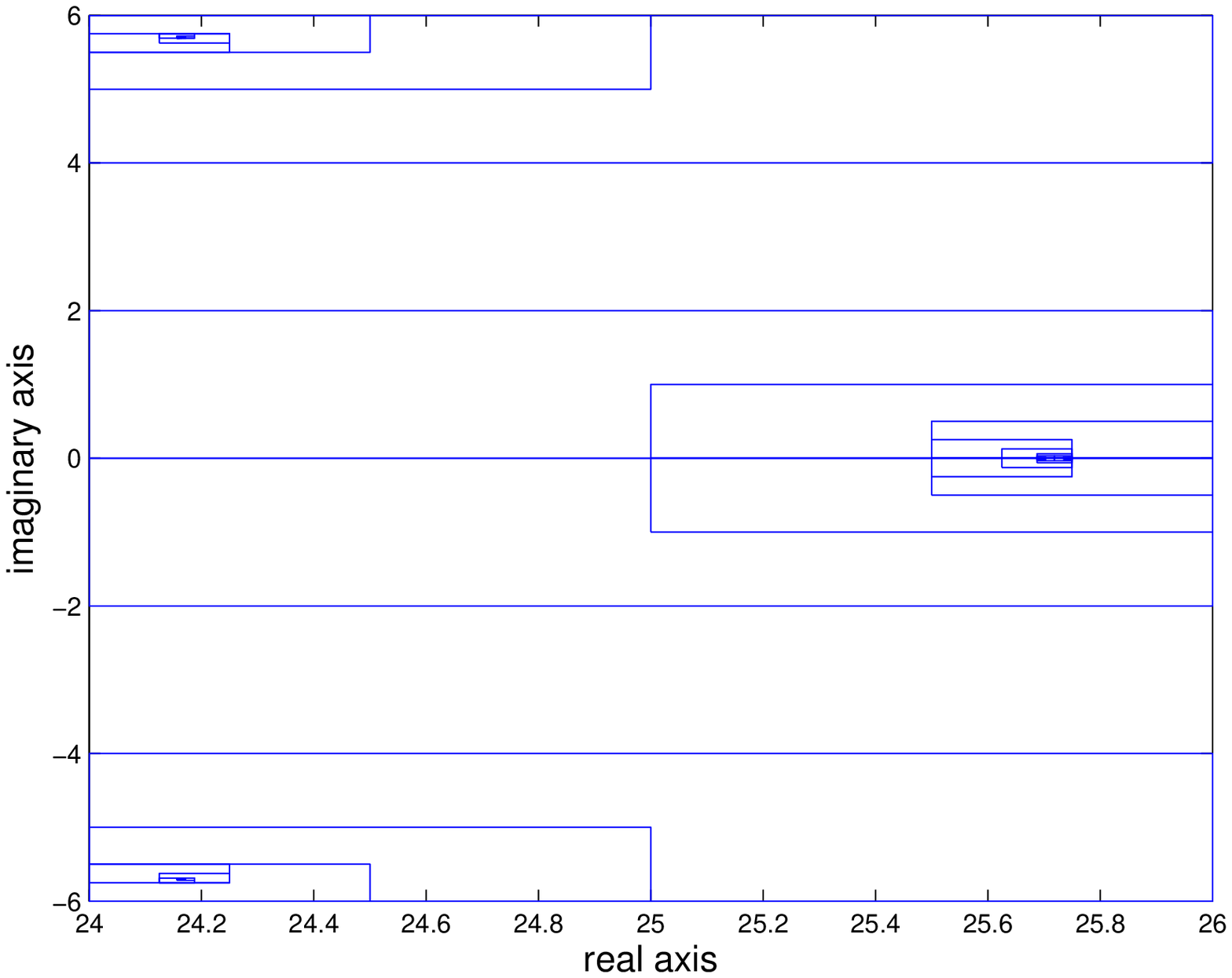}}\\
\end{tabular}
\end{center}
\caption{The regions explored by {\bf RIM} for the disc with radius $1/2$, $n(x)=16$, and $\epsilon = 1.0e-3$. Left: the search region is given by $S=[1, 10] \times [-1, 1]$. Right: the search region is given by $S=[22, 25] \times[-8, 8]$.}
\label{circlen16}
\end{figure}

The search regions for the transmission eigenvalue tests 
are shown in Fig \ref{circlen16}. The algorithm refines 
near the eigenvalues until the tolerance is met.
The right image in Fig. \ref{circlen16} shows 
only three refined regions because two eigenvalues are 
very close.

{\bf Example 2:} Let $D$ be the unit square and $n(x)=16$
(AAS: I think this is correct) with $h\approx0.05$. 
The matrices $A$ and $B$ are
$2075 \times 2075$. The exact transmission 
eigenvalues are not available.
The first search region is given by $S=[3, 8] \times[-1, 1]$. {\bf RIM} computes the following eigenvalues
\[
\lambda_1 = 3.561, \quad \lambda_2 = 6.049, \quad \lambda_3 = 6.051.
\]
They are consistent with the values given in Table 3 of \cite{CakoniEtal2010IP}:
\[
\lambda_1 = 3.479, \quad \lambda_2 = 5.883, \quad \lambda_3 = 5.891.
\]
The second search region is given by $S=[20, 25] \times [-8, 8]$. The eigenvalues we obtain are 
\[
\lambda_1 = 20.574+5.128 i, \quad \lambda_2 =  20.574-5.128 i, \quad \lambda_3 = 21.595, \quad \lambda_4 = 23.412.
\]

\begin{figure}[hh]
\begin{center}
\begin{tabular}{lll}
\scalebox{0.42}{\includegraphics{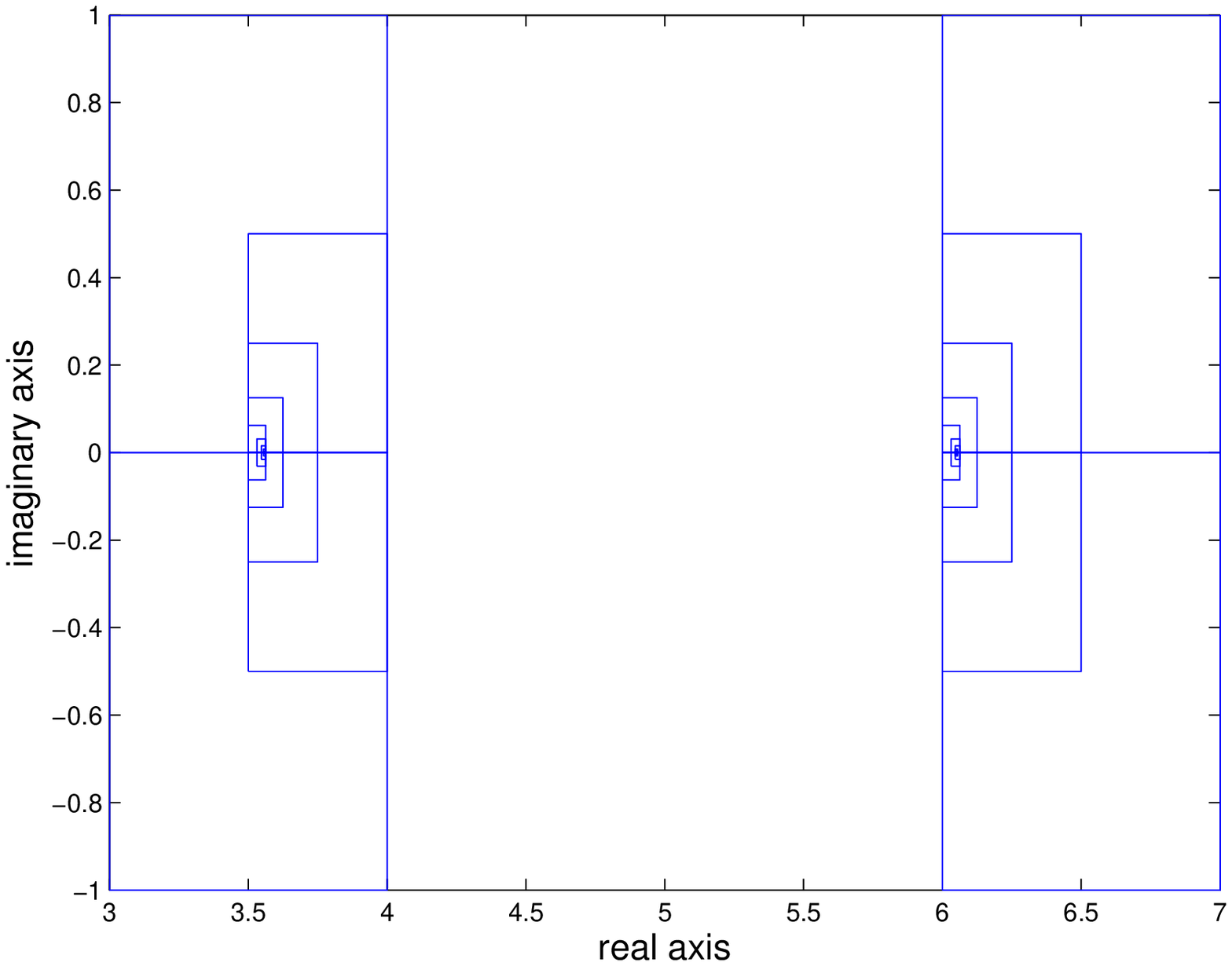}}&
\scalebox{0.42}{\includegraphics{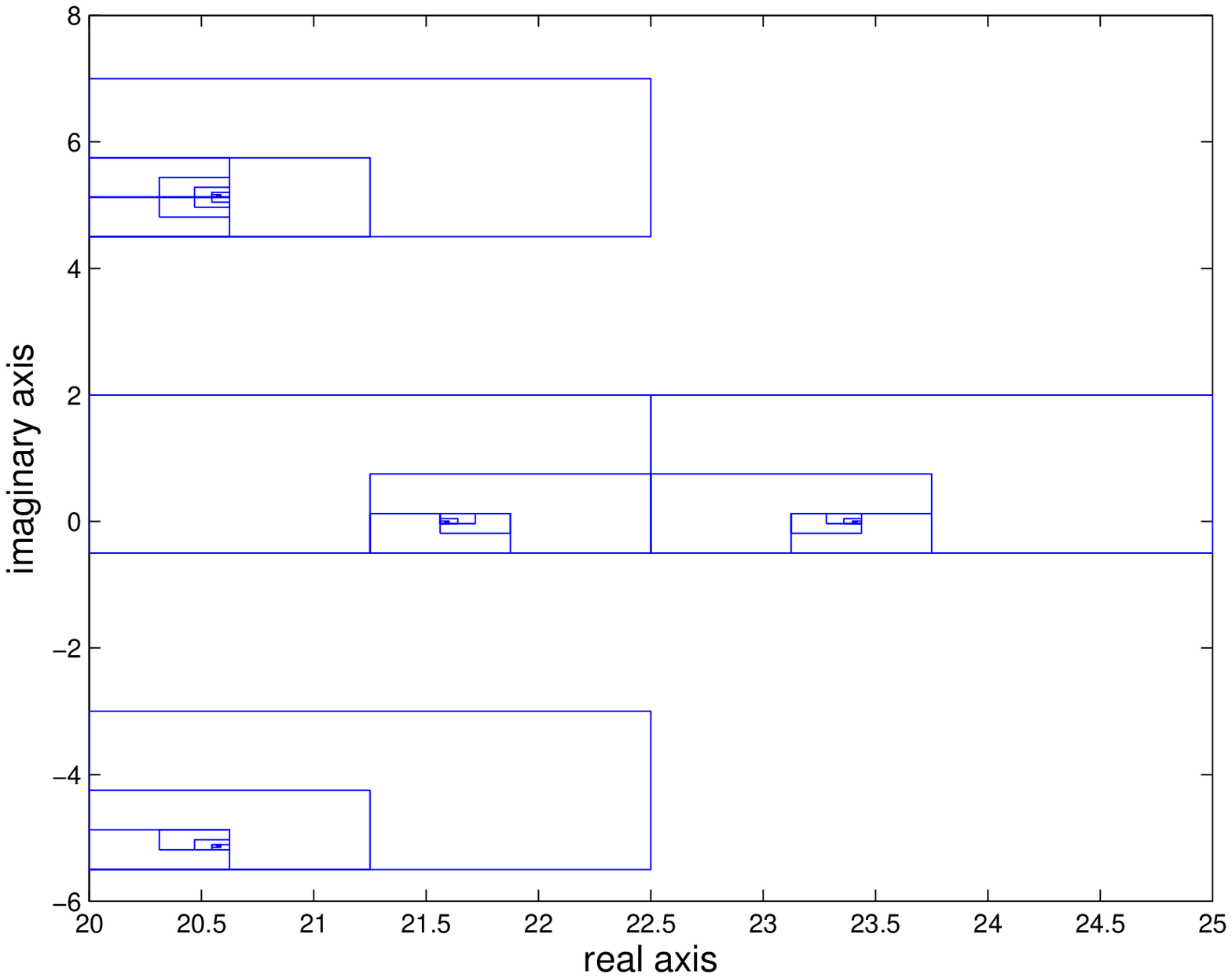}}\\
\end{tabular}
\end{center}
\caption{The regions explored by {\bf RIM} for the unit square with $n(x)=16$ and $\epsilon = 1.0e-3$. 
Left: the search region is given by $S=[3, 7] \times [-1, 1]$. Right: the search region is given by $[20, 25] \times [-6, 8]$.}
\label{square16}
\end{figure}

We plot the search regions in Fig.~\ref{square16}. The left picture is for $S=[3, 7] \times[-1, 1]$. The right picture is for $S=[20, 25] \times [-6, 8]$.

\subsection{Eigenvalues on $\Gamma :=\partial S$}
It is very unlikely that $\Gamma := \partial S$ is not 
contained in the resolvent set. However, we want to 
explore what will happen if eigenvalues lie on on $\Gamma$.
The first example shows that this does not generate difficulty for {\bf RIM}.

{\bf Example 3:} We first consider a simple example given below (Example 5 of \cite{SakuraiSugiura2003CAM}):
\[
A = \begin{pmatrix} \frac{99}{100}& \frac{1}{100}&0&\ldots & 0 \\
				0& \frac{98}{100}&0&\ldots& 0 \\
				\ddots & \ddots& \ddots &\ddots &\vdots \\
				0&\ldots&0&\frac{1}{100}&\frac{1}{100} \\
				0&\ldots&\ldots&0&\frac{0}{100}
				\end{pmatrix}, \quad
B = \text{diag}(0,\ldots,0, 1, \ldots, 1),
\]
where $B$ has $20$ ones on its diagonal. The following are some exact eigenvalues
\[
\lambda_1 = 0, \quad \lambda_2=0.01,\quad\lambda_3=0.02, \quad \lambda_4=0.03.
\]

We set the initial search region to be 
$S=[0, 1/30] \times[0, 1/100]$ and
$\epsilon = 1.0e-9$ and note that all the 
eigenvalues are on $\Gamma:=\partial S$.

\begin{figure}
\begin{center}
{ \scalebox{0.6} {\includegraphics{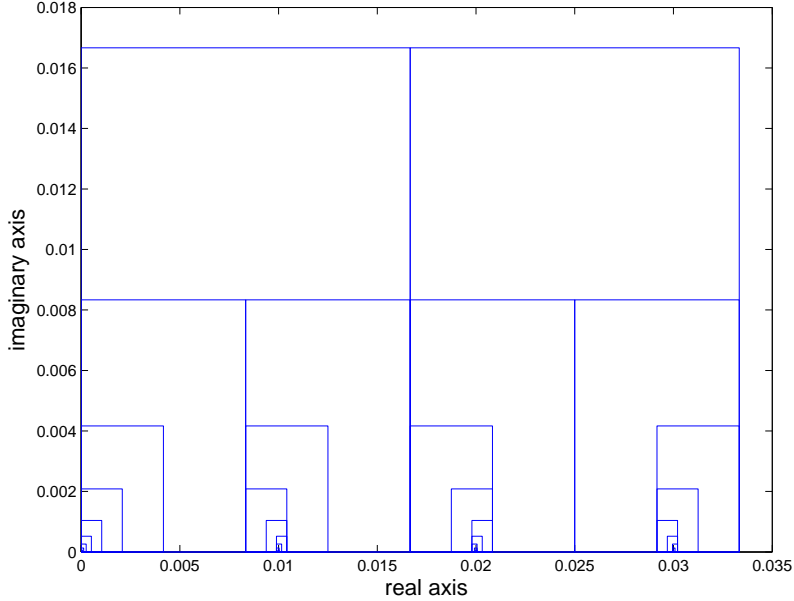}}}
\caption{Eigenvalues on $\Gamma=\partial S$. All the four eigenvalues are on $\Gamma$.}
\label{SS5}
\end{center}
\end{figure}

The eigenvalues computed by {\bf RIM} are
given below (see also Fig. \ref{SS5}). They are accurate up to the required precision. From Fig. \ref{SS5}, we can see that {\bf RIM}
keeps refining around the eigenvalues.
\begin{eqnarray*}
&&\lambda_1 = (4.967053731282552 + 4.967053731282552i)10^{-10},\\ 
&&\lambda_2= 0.009999999900659 + 0.000000000496705i, \\
&& \lambda_3=0.020000000298023 + 0.000000000496705i, \\
&& \lambda_4 =0.029999999701977 + 0.000000000496705i.
\end{eqnarray*}

\subsection{Self-correction Property}
When a quadrature point $z_q$ in the collection
of linear systems (\ref{projectionr}) is close 
to an eigenvalue $\lambda$, the linear system will 
be ill-conditioned. In particular, when $\lambda$ is 
just outside $\Gamma$ the indicator function 
$\chi_S$ could be large because either the linear solve
or quadrature rule are not sufficiently accurate. 
{\bf RIM} will take such regions as admissible and
refine. But fortunately, after a few subdivisions, 
{\bf RIM} appears to discard the sub regions. We 
demonstrate this interesting {\it self-correction 
property} using two example.

{\bf Example 4:} We use matrices $A$ and $B$ 
from Example 2 and focus on the
eigenvalue located at $3.9945$. 
We choose the initial search region 
$S = [4.0, 4.2] \times [0, 0.2]$ and note that 
there is no eigenvalue in $S$. With the same
standard two-point Gauss quadrature rule 
on each edge of $S$ {\bf RIM} computes
\begin{equation}\label{chiS}
\chi_{S} = 4.383,
\end{equation}
indicating that there may be eigenvalues in $S$
and {\bf RIM} procedes to recursively explore
$S$ by dividing $S$ into the four rectangles
\begin{eqnarray*}
&&S^1_1=[4.0, 4.1] \times [-0.1, 0], \quad S^1_2=[4.0, 4.1] \times [0, 0.1],  \\
&&S^1_3=[4.1, 4.2] \times [-0.1, 0], \quad S^1_4=[4.1, 4.2] \times [0, 0.1]
\end{eqnarray*}
with indicator values
\[
\chi_{S^1_1} = 1.589, \quad \chi_{S^1_2} = 1.589, \quad \chi_{S^1_3}= 0.002, \quad \chi_{S^1_4} = 0.002.
\]
{\bf RIM} discards $S^1_3$ and $S^1_4$ and retains 
$S^1_1$ and $S^1_2$ as admissible regions.

We show the result for region $S^1_1$: $S^1_2$ is similar. 
The four rectangles from dividing $S_1^1$ are
\begin{eqnarray*}
&&S^2_1 = [4.0, 4.05] \times [-0.05, 0], \quad S^2_2 = [4.0, 4.05] \times [-0.1, -0.05], \\
&&S^2_3 = [4.05, 4.1] \times [-0.05, 0], \quad S^2_4 = [4.0, 4.05] \times [-0.1, -0.05],
\end{eqnarray*}
with indicator values
\[
\chi_{S^2_1} = 0.997, \quad \chi_{S^2_2} = 0.002, \quad \chi_{S^2_3}= 0.002, \quad \chi_{S^2_4} = 2.159e-04.
\]
and {\bf RIM} discards all the regions. Let us see one more level. Suppose
$\chi_{S^2_1}$ is kept and subdivided into
\begin{eqnarray*}
S^3_1 = [4.0, 4.025] \times [-0.025, 0], \quad S^3_2 = [4.0, 4.025] \times [-0.05, -0.025], \\
S^3_3 = [4.025, 4.05] \times [-0.025, 0], \quad S^3_4 = [4.025, 4.05] \times
[-0.05, -0.025]
\end{eqnarray*}
with indicator values
\[
\chi_{S^3_1} = 0.395, \quad \chi_{S^3_2} = 0.002, \quad \chi_{S^3_3}= 0.001,
\quad \chi_{S^3_4} = 1.615e-04.
\]
Hence {\bf RIM} eventually discards $S$.

{\bf Example 5:} The same experiment is conducted
for a search region around the complex
eigenvalue $\lambda = 24.1586 + 5.690i$
with initial search region 
$S=[24.16, 24.96] \times [5.30, 6.10]$
which although close to the eigenvalue
does not contain any eigenvalues.
Indicator values are in Table.
\ref{selfcorrectionA} and we can note that {\bf RIM} does
eventually conclude that there are no eigenvalues
in the region.
\begin{center}
\begin{table}[hh]
\caption{The indicators function $\chi_S$ on different search regions. }
\label{selfcorrectionA}
\begin{center}
\begin{tabular}{lrlr}
\hline
$S^1_1=[24.16, 24.56] \times  [5.30, 5.70]$ &11.825& $S^1_2=[24.16, 24.56] \times  [5.70, 6.10]$ &0.195 \\
$S^1_3=[24.56, 24.96] \times  [5.30, 5.70]$ &5.418e-11& $S^1_4=[24.56, 24.96] \times  [5.70, 6.10]$ &4.119e-11 \\
\hline
$S^2_1=[24.16, 24.36] \times  [5.30, 5.50]$ & 9.216e-11& $S^2_2=[24.16, 24.36] \times  [5.50, 5.70]$&3.682 \\
$S^2_3=[24.36, 24.56] \times  [5.30, 5.50]$& 8.712e-14& $S^2_4=[24.36, 24.56] \times  [5.50, 5.70]$&5.870e-11 \\
\hline
$S^3_1=[24.16, 24.26] \times  [5.50, 5.60]$ &1.742e-11& $S^3_2=[24.16, 24.26] \times  [5.60, 5.70]$&7.806 \\
$S^3_3=[24.26, 24.36] \times  [5.50, 5.60]$ &1.476e-13 & $S^3_4=[24.26, 24.36] \times  [5.60, 5.70]$ &6.755e-11 \\
\hline
$S^4_1=[24.16, 24.21] \times  [5.60, 5.65]$ &6.558e-10& $S^4_2=[24.16, 24.21] \times  [5.65, 5.70]$ &2.799 \\
$S^4_3=[24.21, 24.26] \times  [5.60, 5.65]$ &1.378e-13& $S^4_4=[24.21, 24.26] \times  [5.65, 5.70]$ &8.229e-11 \\
\hline
$S^5_1=[24.16, 24.185] \times  [5.65, 5.675]$ &1.159e-8& $S^5_2=[24.16, 24.185] \times  [5.675, 5.70]$ &1.556\\
$S^5_3=[24.185, 24.21] \times  [5.65, 5.675]$ &4.000e-13& $S^5_4=[24.185, 24.21] \times  [5.675, 5.70]$ &8.648e-11\\
\hline
$S^6_1=[24.16, 24.185] \times  [5.65, 5.675]$ &5.574e-06& $S^6_2=[24.16, 24.1725] \times  [5.6875, 5.70]$ &0.095\\
$S^6_3=[24.185, 24.21] \times  [5.65, 5.675]$ &4.304e-12& $S^6_4=[24.185, 24.21] \times  [5.675, 5.70]$ &2.628e-11\\
\hline
\end{tabular}
\end{center}
\end{table}
\end{center}

\subsection{Close Eigenvalues}
{\bf RIM} is able to separate nearby eigenvalues 
provided the tolerance is less than the eigenvalue 
separation.

{\bf Example 6:} This example comes from a 
finite element discretization of the
Neumann eigenvalue problem:
\begin{subequations}\label{Neumann}
\begin{align}
 - \triangle u&=\lambda u,  &\text{in } D,\\[1mm]
 \frac{\partial u}{\partial \nu} &=0,  &\text{on } \partial D,
\end{align}
\end{subequations}
where $D$ is the unit square which has an eigenvalue $\pi^2$ 
of multiplicity $2$. We use linear Lagrange elements 
on a triangular mesh with $h\approx 0.025$
to discretize and obtain a generalized eigenvalue problem
\begin{equation}\label{FEMNeumann}
A{\bf x} = \lambda B {\bf x},
\end{equation}
where the stiffness matrix $A$ and
mass matrix $B$ are $2075 \times 2075$.
The discretization has broken the symmetry and
\eqref{FEMNeumann} the eigenvalue of multiplicity
$2$ has been approximated by a very close
pair of eigenvalues of
\[
\lambda_1=9.872899741642826 \quad \text{and} \quad \lambda_2=9.872783160389966.
\] 

With  $\epsilon = 1.0e-3$ ${\bf RIM}$ fails to 
separate the eigenvalues and we obtain only one 
eigenvalue
\[
\lambda_1 = 9.872680664062500.
\]
However, with $\epsilon = 1.0e-9$
${\bf RIM}$ separates the eigenvalues and we obtain 
\[
\lambda_1=9.872899741516449 \quad \text{and} \quad \lambda_2 = 9.872783160419203.
\]
The search regions explored by {\bf RIM} with different tolerances are shown in Fig.~\ref{Neumann}.
\begin{figure}[hh]
\begin{center}
\begin{tabular}{lll}
\scalebox{0.42}{\includegraphics{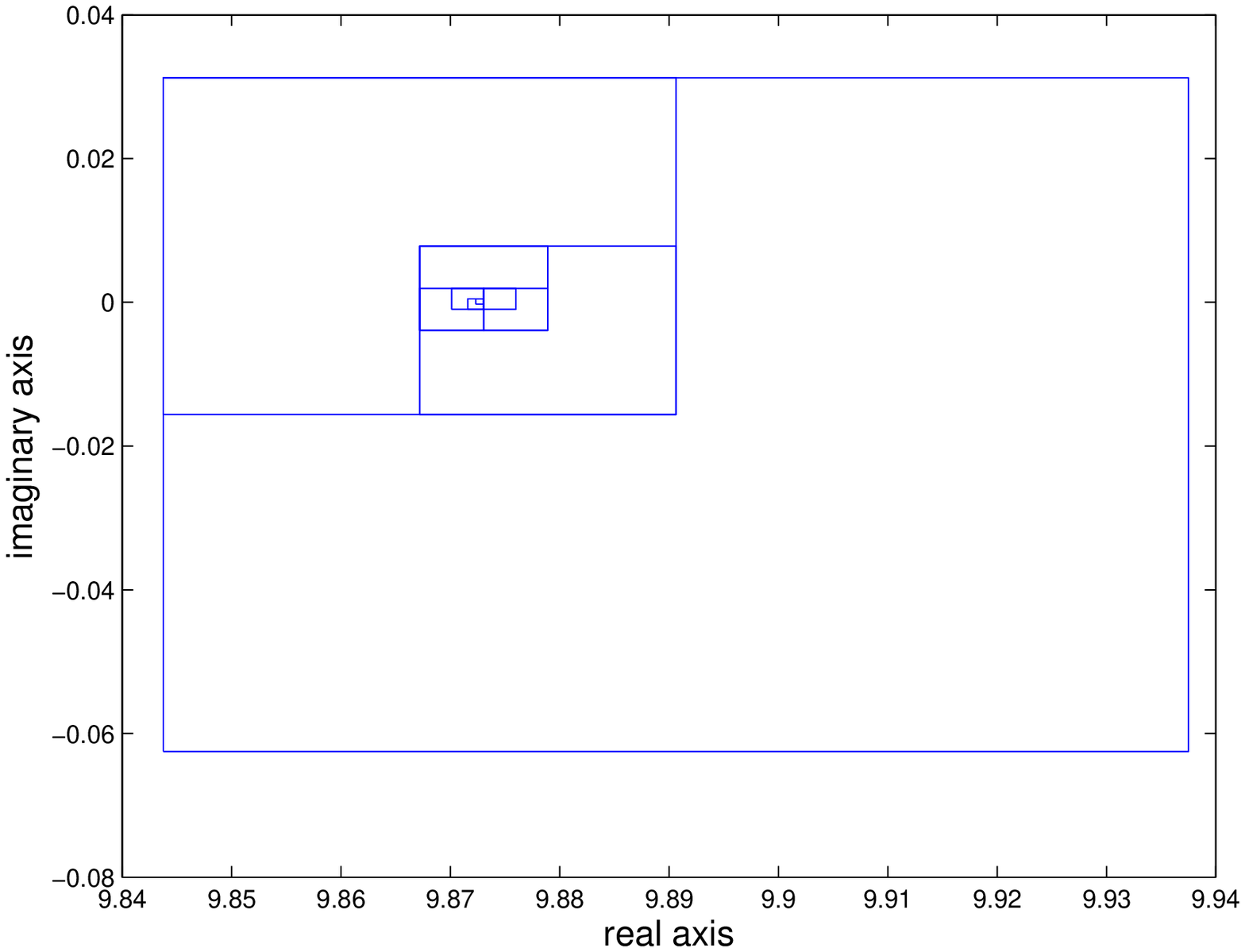}}&
\scalebox{0.42}{\includegraphics{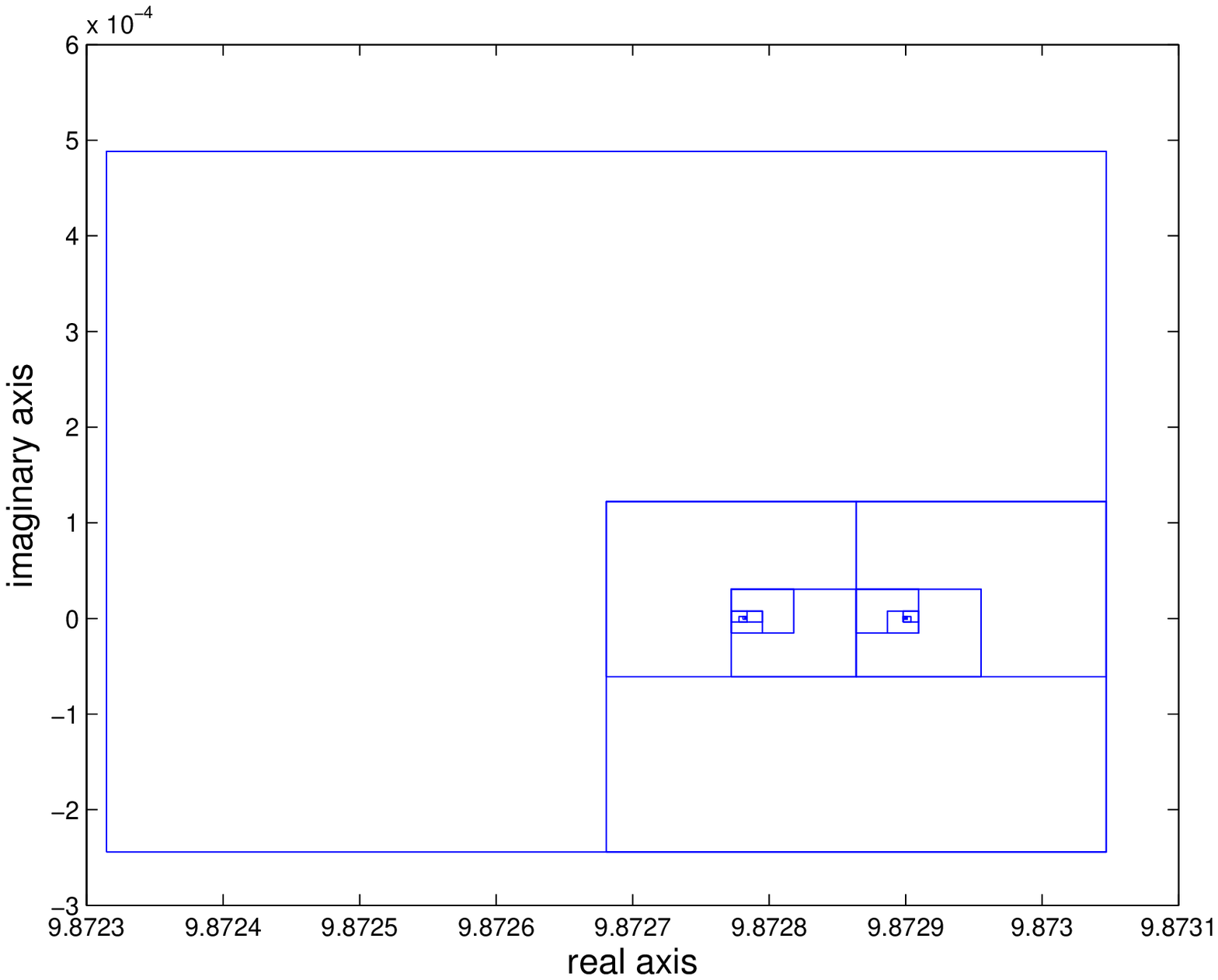}}\\
\end{tabular}
\end{center}
\caption{The regions explored by {\bf RIM}. The search region is given by $S=[1, 10] \times [-1, 1]$. Left: $\epsilon =1.0e-3$. 
Right: $\epsilon = 1.0e-9$.}
\label{Neumann}
\end{figure}

\begin{figure}
\begin{center}
{ \scalebox{0.6} {\includegraphics{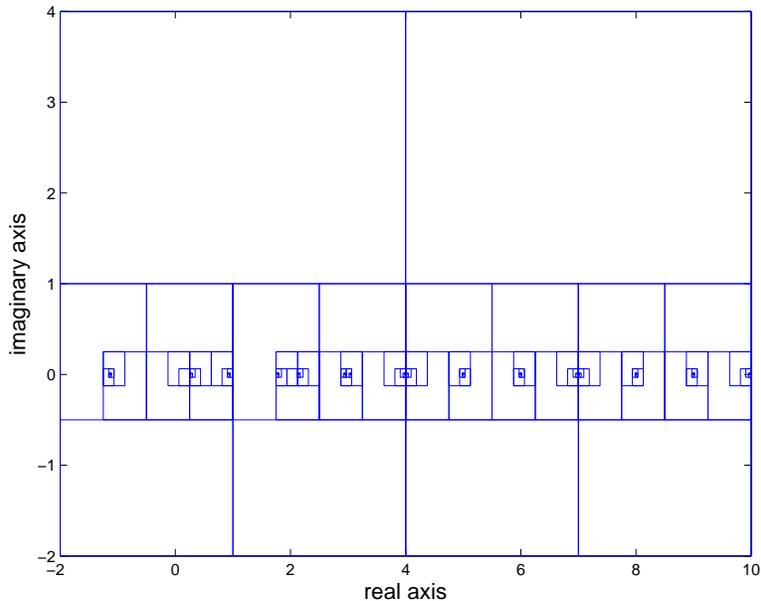}}}
\caption{The regions explored by {\bf RIM} for the Wilkinson matrix with $\epsilon = 1.0e-14$.}
\label{Wilkinson}
\end{center}
\end{figure}

{\bf Example 7:} As a final example we compute 
the eigenvalues of the $40 \times 40$ Wilkinson matrix
\[
A = \begin{pmatrix} 19& -1 \\ 
				-1& 18 &-1 \\
				\quad& \ddots &\ddots & \ddots \\ 
				 \quad&\quad&-1&1&-1  \\  
				 \quad&\quad&\quad&-1&1&-1  \\
	 			\quad&\quad&\quad&\quad&\ddots &\ddots & \ddots \\ 
				 \quad&\quad&\quad&\quad&\quad&-1& 19 &-1 \\ 
				 \quad&\quad&\quad&\quad&\quad&\quad&-1&20\end{pmatrix}
\]
which is known to have very close eigenvalues. 
With 
$\epsilon = 1.0e-14$ and the search region $S = [-2, 10] \times [-2, 10]$.
{\bf RIM} accurately distinguishes the close eigenvalues with
giving the results shown in Table \ref{WilkinsonEig} and Fig. \ref{Wilkinson}.
\begin{center}
\begin{table}[hh]
\caption{The computed Wilkinson eigenvalues by {\bf RIM}.}
\label{WilkinsonEig}
\begin{center}
\begin{tabular}{rlrl}
$i$ & $\lambda_i$ &$i$ & $\lambda_i$\\
\hline
1& -1.125441522046458 &11 &5.000236265619321 \\
2& 0.253805817279499 &12 &5.999991841327017 \\
3&0.947534367500339 &13 &6.000008352188331 \\
4&1.789321352320258 &14 &6.999999794929806 \\
5&2.130209219467361 &15 &7.000000207904748 \\
6&2.961058880959172 &16 & 7.999999996191775\\
7&3.043099288071971 &17 &8.000000003841876\\
8&3.996047997334983 &18 &8.999999999945373 \\
9&4.004353817323874 &19 &9.000000000054399 \\
10&4.999774319815003 &20 &9.999999999999261 \\
\hline
\end{tabular}
\end{center}
\end{table}
\end{center}

\section{Discussion and future work}
This paper proposes a robust recursive integral method 
{\bf RIM} to compute transmission eigenvalues.
The method effectively locates all eigenvalues
in a region when neither the location or number 
eigenvalues is known.  The key difference between 
{\bf RIM} and other counter integral based methods 
in the literature is that {\bf RIM} essentially
only tests if a region contains eigenvalues or not.
As a result accuracy requirements on quadrature,
linear solves, and the number of test vectors
may be significantly reduced. 

{\bf RIM} is a non-classical eigenvalue solver which
is well suited to problems that only require eigenvalues.
In particular, the method snot only works for 
matrix eigenvalue problems resulting from 
suitable numerical approximations,
e.g., finite element methods, of PDE-based
eigenvalue problem, but also those eigenvalue 
problems which can not be easily casted as a 
matrix eigenvalue problem, e.g.,
see \cite{Beyn2012LAA, Kleefeld2013IP}.

The goal of this paper is to introduce the 
idea of {\bf RIM} and demonstrate its 
potential to compute eigenvalues. 
A paper like this raises more questions than 
it answers. How {\it inaccurate} 
can the quadrature be and still locate eigenvalues?
How {\it inaccurate} can the the linear solver can and
still locate eigenvalues. The current implementation
uses a combination of inaccurate quadrature and 
accurate solver: two point Gaussian quadrature on the edges of 
rectangles and the Matlab ``$\backslash$" operator.
These two separate issues can be combined into one 
question: how {\it accurate} does the overall 
procedure have to be to accurately 
distinguish admissible regions. 
These crucial complexity issues are not addressed 
in this current paper.

The example problems are small. We plan to extend 
{\bf RIM} for large (sparse) eigenvalue problems
which will require replacing ``$\backslash$"
with an iterative solver. 
Parallel extension is another interesting 
project since the algorithm is essentially 
{\it embarrassingly parallel}. In
particular, a GPU implement of {\bf RIM}
is under consideration.

\section*{Acknowlegement}
The work of JS and RZ is partially supported NSF DMS-1016092/1321391.

\end{document}